\documentclass[11pt]{article}

\usepackage{amsmath,amsfonts,amssymb,amsthm}
\usepackage[english]{babel}

\newcommand{\sA}{{\mathcal A}}
\newcommand{\sB}{{\mathcal B}}
\newcommand{\sD}{{\mathcal D}}

\newcommand{\sH}{{\mathcal H}}
\newcommand{\sL}{{\mathcal L}}
\newcommand{\sR}{{\mathcal R}}

\newtheorem{theorem}{Theorem}[section]

\newtheorem{lemma}[theorem]{Lemma}

\newtheorem{remark}[theorem]{Remark}

\newtheorem{definition}[theorem]{Definition}

\newtheorem*{problem}{Problem}

\numberwithin{equation}{section}

\begin{document}


\centerline{{\Large \bf Regularity of boundary traces for a fluid-solid \vspace{2mm}}} 
\centerline{{\Large \bf interaction model}\footnote{The research of the first author is supported by the Italian MIUR, within the project 2007WECYEA (``Metodi di viscosit\`a, metrici e di teoria del controllo in equazioni alle derivate parziali non lineari'').
The research of the second author is supported by the NSF grant DMS-060666882 and by the 
AFOSR grant FA9550-09-1-045.  }}

\vspace{5mm}


\centerline{
{\large {\sc Francesca Bucci}$^{\textrm{a}}$, {\sc Irena Lasiecka}$^{\textrm{b}}$}}
\vspace{8mm}      

\noindent
{\small $^{\textrm{a}}$Universit\`a degli Studi di Firenze, Dipartimento di
Matematica Applicata, Via S.~Marta 3, 50139 Firenze, Italy; e-mail
{\tt francesca.bucci@unifi.it}} 
\vspace{2mm}\\  
{\small $^{\textrm{b}}$University of Virginia, Department of
Mathematics, P.O.~Box 400137, Charlottesville, VA 22904-4137, U.S.A.,
e-mail {\tt il2v@virginia.edu}} \vspace{8mm}
  
\noindent
{\small {\bf Abstract:}
We consider a mathematical model for the interactions of an elastic body fully immersed 
in a viscous, incompressible fluid. 
The corresponding composite PDE system comprises a linearized Navier-Stokes system
and a dynamic system of elasticity; the coupling takes place on the interface between the two regions occupied by the fluid and the solid, respectively.   
We specifically study the regularity of boundary traces (on the interface)
for the fluid velocity field.
The obtained trace regularity theory for the fluid component of the system---of 
interest in its own right---establishes, in addition, solvability of the associated optimal (quadratic) control problems on a finite time interval, along with well-posedness of the corresponding operator Differential Riccati equations.
These results complement the recent advances in the PDE analysis and control of the 
Stokes-Lam\'e system. 
\vspace{2mm}
}

\section{Introduction}
The description of the interactions between a fluid flow and an elastic structure via an interface is a classical topic in continuum mechanics.
Motivated by an impressive variety of fields of application, such as, e.g., cell biology, biomedical engineering, ultrasound imaging, in the past decade fluid-structure interaction models have received a renewed attention.  
Numerical and experimental studies have pursued the derivation of appropriate PDE models 
for specific physical problems, along with the understanding of the corresponding dynamics; 
see, e.g., \cite{quarteroni-formaggia} and \cite{moubachir-zol}.
On the other hand, basic mathematical questions such as well-posedness (in a natural functional setting) of the systems of coupled partial differential equations (PDE) which result from the modeling of these interactions, have been explored and answered only recently.

Let us recall that a major technical issue which arises in the PDE analysis of 
fluid-solid interaction models is the apparent discrepancy between the regularity of boundary traces for the fluid and the elastic (hyperbolic) components on the interface, 
with the latter not defined {\em a priori} by means of standard trace theory.
Different solutions have been proposed to cope with the mismatch between parabolic and hyperbolic regularity. For instance, in \cite{boulakia} a regularizing term is included in the elastic equation, thereby ensuring the soughtafter boundary regularity, as a consequence of the improved interior regularity; the works \cite{coutand-shkoller} and \cite{feireisl} consider instead very smooth data, thus obtaining only {\em local} (in time) existence.
In contrast with previous work, the study of \cite{barbu-etal-1} achieves---at least in 
the case of static interface---well-posedness of an established nonlinear fluid-structure interaction model in the natural energy space, along with the exceptional yet intrinsic regularity of the boundary traces of the (elastic) hyperbolic component of the system on the interface.
(To the authors' best knowledge, a similar well-posedness result for the more 
challenging case of moving interface is, at present, lacking.)
  
The composite PDE system under consideration, which is problem 
\eqref{e:navierstokes-0} in the next section, is a linearized version 
of the nonlinear system studied in \cite{barbu-etal-1,barbu-etal-2}; from a 
modeling point of view, the linear evolution is consistent with the 
case of very {\em slow} motion of the fluid.
The boundary value problem comprises a linearized Navier-Stokes system (often referred to 
as Oseen equations) for the velocity field $u$ of the fluid flow and a dynamic system of elasticity for the displacement $w$ of the solid.  
The interaction between the two means takes place at an interface between the two regions occupied by the fluid and the solid, respectively.
We consider the case in which the motion of the solid is entirely due to infinitesimal displacements, and therefore the fluid-solid interface, later denoted by $\Gamma_s$, is  stationary. 
(Yet, although the displacement $u$ is small, the velocity $u_t$ is not,
resulting in the interface condition $u=w_t$.)
While the PDE problem \eqref{e:navierstokes-0} is present in the literature since the 
late sixties (see \cite{lions}), its well-posedness was not shown until 
\cite{barbu-etal-1}, whose theory encompasses the present model.
In the case of the (linear) Stokes-Lam\'e system, the contemporary work of 
\cite{avalos-trig-1,avalos-trig-overview,avalos-trig-3} provides not only well-posedness, 
but also a careful stability analysis.
A very nice overview of recent (theoretical as well as numerical) contributions to 
this topic is found in \cite{du-etal}, whose introduction also provides insight into modeling aspects; see the references therein.

The present study focuses on the regularity of boundary traces for the 
velocity field of the fluid flow $u$ (and $u_t$) on the interface $\Gamma_s$.
Although the obtained regularity properties are of interest in their own right, 
they are also critical in the study of the optimal boundary control problems 
(with quadratic functionals) associated with the fluid-solid interaction 
\eqref{e:navierstokes-0}, in the presence of active controls on the interface 
$\Gamma_s$, thus constituting its main motivation as well as application.

It was shown in \cite{las-tuff-3} that the fluid-structure interaction model under 
consideration yields a suitable (`singular') estimate for the corresponding 
abstract evolution, which plays a major role in the study of the associated 
optimal control problems.
Indeed, in light of the regularity result obtained in \cite[Theorem~5.1]{las-tuff-3}, the optimal control theory in \cite{las-cbms,las-trig-se-2,las-tuff-2} applies, ensuring a feedback control law with a {\em bounded} gain operator (defined on the state space), along with well-posed Riccati equations. 
This, however, holds true assuming that the functional penalizes the quadratic energy 
of the solid component in a weaker topology than the actual energy level (this 
is a consequence of the fact that the singular estimate holds in the norm
of a function space $\sH_{-\alpha}$ which is larger than the energy space). 

A first attempt to remove the constraint on the structure of the quadratic cost 
functional---which corresponds to a constraint on the observation operator---has 
been successfully carried out in \cite{bucci-las-cvpde} in the case of the simpler 
Stokes-Lam\'e system.
The PDE analysis performed in \cite{bucci-las-cvpde} shows indeed that the singular 
estimate in \cite[Theorem~5.1]{las-tuff-3} does not hold for $\alpha=0$ and hence 
that the optimal boundary control theory of \cite{las-cbms,las-tuff-2} is
inapplicable to the case of natural cost functionals such as the integral of the 
quadratic energy of the system.
Nonetheless, \cite{bucci-las-cvpde} provides novel regularity results which 
ensure solvability of the associated optimal control problems, with well-posed Riccati equations, according to the distinct (optimal control) theory contained in \cite{abl-2}.
 
In the present paper we show that the complex of boundary regularity results 
established in \cite{bucci-las-cvpde} for the uncontrolled Stokes-Lam\'e system 
may be extended to more general fluid-structure interaction models 
comprising a linearized Navier-Stokes (rather than the Stokes) system for a slow, 
viscous, incompressible fluid (see \eqref{e:navierstokes-0} below).
As we shall see, the regularity of boundary traces for the fluid component of the system 
will follow combining the exceptional boundary regularity of the hyperbolic component of 
the system ({\em cf}.~ Lemma~\ref{l:wave-regularity}), the regularizing effect of the 
parabolic component, and interpolation techniques.
The same factors are behind the proof of the specific regularity estimates recently derived for diverse composite PDE systems, including thermoelastic systems 
(\cite{bucci-las-thermo,abl-1}) and acoustic-structure interaction models 
(\cite{bucci-applicationes}).

\smallskip
An outline of the paper follows below.
In Section~\ref{s:setting-statements} we introduce the coupled PDE system, along with 
the natural functional setting where well-posedness holds true, in both the variational 
and semigroup sense. 
The statement of the present work's main contribution, namely Theorem~\ref{t:main}, 
which collects several boundary regularity results for the fluid velocity field of 
problem \eqref{e:navierstokes-0}, is given here.
Section~\ref{s:known-results} contains some preliminary trace regularity results which are 
crucially instrumental to the proof of Theorem~\ref{t:main}. 
Section~\ref{s:proof} is entirely devoted to the proof of all the assertions of Theorem~\ref{t:main}.
Application of the obtained regularity theory to the associated optimal control problems is discussed in Section~\ref{s:application}, which concludes the paper.

\section{The mathematical model, statement of the main result}
\label{s:setting-statements}
In this section we introduce the PDE system under investigation.
We preliminarly recall from \cite{barbu-etal-1} the well-posedeness and regularity 
theory which constitutes the basis of our analysis.
Theorem~\ref{t:main} contains the statement of the novel regularity results established 
for the fluid component of the system.
\subsection{The PDE problem}
Let $\Omega_f$ and $\Omega_s$ be the open smooth (three-dimensional) domains 
occupied by the fluid and the solid, respectively.
The boundary of $\Omega_s$ represents the {\em interface} between the fluid and the solid, and is denoted by $\Gamma_s=\partial\Omega_s$;
$\Gamma_f$ is the outer boundary of $\Omega_f$, that is
$\Gamma_f=\partial\Omega_f\setminus \partial\Omega_s$.
It is assumed that the motion of the solid is entirely due to infinitesimal 
displacements, and hence that the interface $\Gamma_s$ is {\em stationary}.
The PDE system under investigation comprises a linearized Navier-Stokes system for
the fluid velocity $u$ in $\Omega_f$, 
and an elastic equation for the displacement $w$ of the solid in 
$\Omega_s$. 
The coupling takes place on the interface $\Gamma_s$; the case of small, but rapid oscillations of the elastic body is considered, which yields the interface condition 
$u=w_t$ on $\Gamma_s$ in lieu of the no-slip boundary condition $u=0$.
(For a detailed discussion of modeling issues see, e.g. \cite{du-etal}.)  
Thus, if $p$ denotes the fluid pressure, the pair $(u,w,p)$ satisfies
\begin{equation} \label{e:navierstokes-0}
\left\{ \hspace{1mm}
\begin{split} 
& u_t-{\rm div}\,\epsilon(u) + \sL u
+\nabla p= 0 & &\textrm{in }\; Q_f:= \Omega_f\times (0,T)
\\
& {\rm div}\, u=0 & & \textrm{in }\; Q_f
\\
& w_{tt} - {\rm div}\,\sigma(w)=0 & &\textrm{in }\; Q_s:= \Omega_s\times (0,T)
\\
& u=0 & & \textrm{on }\; \Sigma_f:= \Gamma_f\times (0,T)
\\
& w_t=u & & \textrm{on }\; \Sigma_s:= \Gamma_s\times (0,T)
\\
& \sigma(w)\cdot \nu=\epsilon(u)\cdot\nu - p\nu& & 
\textrm{on }\; \Sigma_s
\\
& u(\cdot,0)=u_0 & & \textrm{in }\; \Omega_f
\\
& w(\cdot,0)=w_0\,, \quad w_t(\cdot,0)=z_0 & & \textrm{in }\; \Omega_s\,.
\end{split}
\right.
\end{equation}
In the above PDE model $\epsilon$ and $\sigma$ denote the strain 
tensor and the elastic stress tensor, respectively, that are
\begin{equation}
\epsilon_{ij}(\omega)= \frac12\Big(\frac{\partial \omega_i}{\partial x_j}
+\frac{\partial \omega_j}{\partial x_i}\Big)\,,
\qquad
\sigma_{ij}(\omega) = \lambda \sum_{k=1}^3 \epsilon_{kk}(\omega)\,\delta_{ij}
+ 2\mu\,\epsilon_{ij}(\omega)\,,
\end{equation}
where $\lambda, \mu>0$ are the Lam\'e constants and $\delta_{ij}$ is the Kronecker
symbol.
The (uncoupled) linear system for the fluid velocity field, often referred to 
as Oseen equations, is the linearization of the Navier-Stokes equations about 
a given (time-independent) smooth 
vector valued function $v$ such that $\textrm{div}\,v=0$;
see, e.g., \cite{fursikov-etal}.
Thus, the operator $\sL$ is defined as follows: 
\begin{equation}
\sL u(x,t)\equiv (v(x)\cdot\nabla)u(x,t)+ (u(x,t)\cdot\nabla)v(x)\,.
\end{equation}

\subsection{Functional setting, well-posedness}
The energy space for the PDE problem \eqref{e:navierstokes-0} is 
\begin{equation*}
Y = H\times H^1(\Omega_s)\times L_2(\Omega_s)\,,
\end{equation*}
where $H$ is defined as follows:
\begin{equation*}
H := \big\{ u\in L_2(\Omega_f): \; {\rm div}\, u=0\,, \, u\cdot \nu|_{\Gamma_f}=0\big\}\,,
\end{equation*}
where the divergence is meant in a distributional sense.
$P_H$ will denote the Leray's projection operator from $L_2(\Omega_f)$ onto $H$;
then, $L=P_H\sL$.
In addition, $V$ denotes the space defined as follows:
\begin{equation*}
V := \big\{ u\in H^1(\Omega_f): \; {\rm div}\, u=0\,, \, u|_{\Gamma_f}=0\big\}\,;
\end{equation*} 
The following distinct inner products will be used throughout the paper:
\begin{equation*}
(u,v)_f := \int_{\Omega_f}uv \,{\rm d}\Omega_f\,,\quad 
(u,v)_s := \int_{\Omega_s}uv \,{\rm d}\Omega_s\,,\quad
\langle u,v\rangle := \int_{\Gamma_s} uv \,{\rm d}\Gamma_s\,.
\end{equation*}
The space $V$ is topologized with respect to the inner product given by 
\begin{equation*}
(u,v)_{1,f} := \int_{\Omega_f}\epsilon(u)\epsilon(v) d\Omega_f\,;
\end{equation*}
the corresponding (induced) norm $|\cdot|_{1,f}$ is equivalent to the 
usual $H^1(\Omega_f)$ norm, in view of Korn inequality and the Poincar\'e inequality.
The norm $\|\cdot\|_{H^r(D)}$ in the Sobolev space $H^r(D)$
will be denoted by $|\cdot|_{r,D}$ throughout.

\begin{remark}
\begin{rm} 
Note that all the Sobolev spaces $H^r$ related to $u$ and $w$
are actually $(H^r)^3$: the exponent is omitted just for simplicity
of notation.
\end{rm}
\end{remark}
Let us recall from \cite{barbu-etal-1} the definition of {\em weak} solutions to 
the (uncontrolled) PDE system \eqref{e:navierstokes-0}.

\begin{definition}[Weak solution] 
\label{def:weak-sol}
Let $(u_0,w_0,z_0)\in H\times H^1(\Omega_s)\times L_2(\Omega_s)$ and $T>0$. We say that a triple $(u,w,w_t)\in C([0,T],H\times H^1(\Omega_s)\times L_2(\Omega_s))$
is a weak solution to the PDE system \eqref{e:navierstokes-0} if
\begin{itemize}
\item
$(u(\cdot,0),w(\cdot,0),w_t(\cdot,0))=(u_0,w_0,z_0)$,
\item
$u\in L_2(0,T;V)$,
\item
$\sigma(w)\cdot \nu \in L_2(0,T;H^{-1/2}(\Gamma_s))$,
$\frac{d}{dt}w|_{\Gamma_s}=u|_{\Gamma_s}\in L_2(0,T;H^{1/2}(\Gamma_s))$, and
\item
the following variational system holds a.e. in $t\in (0,T)$:
\begin{equation} \label{e:variational-system}
\begin{cases}
\frac{d}{dt}(u,\phi)_f + (\epsilon(u),\epsilon(\phi))_f
+ (L u,\phi)_f - \langle\sigma(w)\cdot \nu,\phi\rangle=0
\\[1mm]
\frac{d}{dt}(w_{t},\psi)_{s} + (\sigma(w),\epsilon(\psi))_{s}
-\langle\sigma(w)\cdot \nu,\psi\rangle=0\,,
\end{cases}
\end{equation}
for all test functions $\phi\in V$ and $\psi\in H^1(\Omega_s)$.
\end{itemize}
\end{definition}

Existence of weak (global) solutions to the PDE problem \eqref{e:navierstokes-0}
follows by a more general---because pertaining to the actual nonlinear model,
which comprises the Navier-Stokes system---result established in \cite{barbu-etal-1}.

\begin{remark}
\begin{rm} 
It is important to emphasise that the regularity of boundary traces 
of the elastic component required by Definition~\ref{def:weak-sol} 
does not follow from the interior regularity 
$(u,w,w_t)\in C([0,T],H\times H^1(\Omega_s)\times L_2(\Omega_s))$ via standard trace theory.
This is an independent regularity result, valid for any solution corresponding to initial data in the finite energy space; see the statement of Theorem~\ref{t:well-posed} below.
\end{rm}
\end{remark}


\begin{theorem}[Existence of weak solutions, \cite{barbu-etal-1}]
\label{t:well-posed}
For any initial datum $(u_0,w_0,z_0) \in Y$ and any $T>0$, there exists a weak 
solution $(u,w,w_t)$ to system \eqref{e:navierstokes-0}, along with 
the following trace regularity for the elastic component of the system:
\begin{equation*}
\nabla w\big|_{\Gamma_s}\in L_2(0,T;H^{-1/2}(\Gamma_s))\,,
\qquad 
\frac{{\rm d}}{{\rm d}t} w\big|_{\Gamma_s}=
w_t\big|_{\Gamma_s}\in L_2(0,T;H^{1/2}(\Gamma_s))\,.
\end{equation*}
\end{theorem}

\smallskip
It is important to recall from \cite{barbu-etal-1} that by properly defining the overall dynamics operator $\sA$, 
the PDE system \eqref{e:navierstokes-0} admits as well a semigroup representation. 
Preliminarly, one needs to introduce the fluid dynamic operator $A: V\to V'$, that is
\begin{equation}\label{e:fluid-operator}
(Au,\phi) = -(\epsilon(u),\epsilon(\phi))
\qquad \forall \,\phi\in V\,,
\end{equation}
and the (Neumann) map $N:L_2(\Gamma_s)\to H$ defined as follows:
\begin{equation}\label{e:neumann-map}
N:g\mapsto h \Longleftrightarrow
(\epsilon(h),\epsilon(\phi))= \langle g,\phi\rangle \qquad \forall \phi\in V\,.
\end{equation}
Then, by using the operators introduced above, the initial/boundary value 
problem \eqref{e:navierstokes-0} becomes the following abstract differential 
system in the variable $y:=(u,w,w_t)$
\begin{equation}\label{e:abstract-system}
y'=\sA y\,, \qquad y(0)=y_0\,,
\end{equation}
where $y_0:=(u_0,w_0,z_0)\in Y$ and $\sA:\sD(\sA)\subset Y\to Y$ is defined by 
\begin{equation}\label{e:sA}
\sA:=
\left(
\begin{array}{ccc}
A-L & AN \sigma(\cdot)\cdot\nu & 0\\
0 & 0 & I\\ 
0 & \textrm{div}\sigma & 0
\end{array}
\right)\,,
\end{equation}
with domain 
\begin{multline}\label{e:dom-sA}
\sD(\sA) := \big\{\,y=(u,w,z)\in Y: \, u\in V, \; A(u+N \sigma(w)\cdot\nu)-Lu\in H,\;
z\in H^1(\Omega_s),\; \\
\textrm{div}\sigma(w) \in L_2(\Omega_s)\,; \; 
z|_{\Gamma_s}=u|_{\Gamma_s}\,\big\}\,.
\end{multline}
It was shown in \cite[Proposition~3.1]{barbu-etal-1} that the operator 
$\sA:\sD(\sA)\subset Y\to Y$ defined by \eqref{e:sA}--\eqref{e:dom-sA}
is the generator of a $C_0$-semigroup $e^{\sA t}$ of contractions on $Y$.
Consequently, any semigroup solution $y:=(u,w,w_t)$ belongs to $C([0,T],Y)$
when $y_0=(u_0,w_0,z_0)\in Y$, while it belongs to $C([0,T],\sD(\sA))$ if
$y_0\in \sD(\sA)$.
A major achievement of the study performed in \cite{barbu-etal-1} is that
it gives, in addition, that any semigroup solution turns out to be as well a 
weak solution of system \eqref{e:navierstokes-0}. 
In particular, it satisfies the variational formulation \eqref{e:variational-system}, 
along with the intrinsic boundary regularity of the stresses on the interface.
\\
A different semigroup set-up, which does not require---at the level of (semigroup) generalized solutions---the present boundary regularity of the normal stresses, has been developed in \cite{avalos-trig-1}; see also \cite{avalos-trig-overview,avalos-trig-3}.  

\begin{remark}\label{r:fractional-powers}
\begin{rm}   
We finally recall from \cite{las-tuff-3} that the operator $A$ defined by
\eqref{e:fluid-operator} may be considered as acting on $H$ with domain 
$\sD(A):=\{\,u\in V:\, |(\epsilon(u),\epsilon(\phi))|\le C|\phi|_H\,\}$.
Readily $A:\sD(A)\subset H\to H$ is a self-adjoint, negative operator and 
therefore it generates an {\em analytic} semigroup $e^{At}$ on $H$.
Then, the fractional powers of $-A$ are well defined.
A well known result relating $(-A)^\alpha$ and the semigroup generated by $A$ is the pointwise (in time) estimate, true for any $0<\alpha \le 1$: 
\begin{equation} \label{e:analytic-estimate}
\|(-A)^\alpha e^{At}\|_{\sL(H)}\sim \frac{C_T}{t^\alpha}\,, \qquad\quad 0<t\le T\,.
\end{equation}
(To simplify the notation, we shall write $A^\alpha$, rather than $(-A)^\alpha$,
throughout; similarly, later ${\sA}^\alpha$ will mean $(-{\sA})^\alpha$.)
\end{rm}
\end{remark}

\subsection{New results: regularity of boundary traces for the fluid component} 
The main contribution of this work is a set of regularity results for 
the traces of the fluid velocity field $u$ (and $u_t$) on the interface 
between the fluid and the solid regions, collected in Theorem~\ref{t:main} below.
While these regularity estimates are of intrinsic interest, they additionally
play a major role in the study of the quadratic optimal control problems associated with the PDE system under investigation, in the presence of control actions on the interface. 
This fact will be briefly discussed in Section~\ref{s:application}.

\begin{theorem}[Main result: boundary regularity for the fluid component] 
\label{t:main}
Let $y(t)=(u(t),w(t),w_t(t))$ be the solution to the coupled PDE system 
\eqref{e:navierstokes-0} corresponding to an initial state $y_0=(u_0,w_0,z_0)$. 
The following assertions pertaining to the regularity of boundary traces for
the fluid velocity field $u$ (and $u_t$) are valid.
\begin{enumerate}
\item
If $y_0\in Y$, then  
\begin{equation}\label{e:u-reg}
u|_{\Gamma_s}\in L_{4-h}(0,T;L_2(\Gamma_s)) + L_p(0,T;L_2(\Gamma_s)) \qquad 
\forall p\ge 1\,,
\end{equation}
for arbitrarily small $h>0$.
More precisely, $u$ admits a decomposition 
\begin{equation}\label{e:decomposition}
u(t) = u_1(t)+u_2(t)\,,
\end{equation}
where $u_1$ satisfies a pointwise (in time) ``singular estimate'', 
namely there exists a positive constant $C_T$ such that
\begin{equation}\label{e:singular-estimate}
\|u_1(t)\|_{L_2(\Gamma_s)} \le 
\frac{C_T}{t^{1/4+\delta}}\|y_0\|_Y \qquad\quad  \forall y_0\in Y\,, \; 
\forall t\in (0,T]\,,
\end{equation}
with arbitrarily small $\delta>0$, while 
\begin{equation}\label{e:u_2-reg}
u_2|_{\Gamma_s}\in L_p(0,T;L_2(\Gamma_s)) \qquad \forall p\ge 1\,.
\end{equation}
\item
If $y_0\in \sD(\sA^{\epsilon})$ for some $\epsilon\in (0,\frac14)$, then
the regularity of $u_2$ in \eqref{e:u_2-reg} improves as follows:
\begin{equation}\label{e:u_2-reg-improved}
u_2|_{\Gamma_s}\in C([0,T],L_2(\Gamma_s))\,.
\end{equation}
\item 
Let now $y_0\in \sD(\sA^{1-\theta})$, with $\theta\in (0,\frac14)$.
Then, there exists $q\in (1,2)$ such that 
\begin{equation}\label{e:u_t-reg}
u_t|_{\Gamma_s} \in L_q(0,T;L_2(\Gamma_s))
\end{equation}
continuously with respect to $y_0$; 
namely, there exists a constant $C_T$ such that 
\begin{equation}  \label{ineq:trace-estimate}
\|u_t\|_{L_q(0,T;L_2(\Gamma_s))} \le C_T\, \|y_0\|_{\sD(\sA^{1-\theta})}\,.
\end{equation}
The exponent $q$ will depend on $\theta$: more precisely, given $\theta\in (0,\frac14)$, 
we have 
\begin{equation} \label{e:range-for-q}
1<q<\frac{4}{3+4\theta}\,.
\end{equation}

\end{enumerate}
\end{theorem}

\subsection{Regularity of the boundary traces for the elastic component}
\label{s:known-results}
We conclude this section by introducing a key regularity result for the normal stresses 
of the elastic component of the system and a trace estimate pertaining to the fluid velocity field (Lemma~\ref{l:wave-regularity} and Lemma~\ref{l:key} below, respectively), 
which are central to the proof of our main result.
We fashion the proof of the assertions in Lemma~\ref{l:wave-regularity} in some details, following the arguments formerly used in \cite[Lemma 5.2]{las-tuff-3}.
Lemma~\ref{l:key} will follow accordingly.
 
\smallskip
The first result stated below is a slight variant of the one given in 
\cite[Lemma 5.2]{las-tuff-3}. 
More precisely, while the regularity established in \cite[Lemma 5.2]{las-tuff-3} 
pertains to the {\em uncoupled} wave equation and is expressed in terms of Sobolev 
norms, the result established in Lemma~\ref{l:wave-regularity} concerns the regularity 
of the elastic component of the coupled system driven by the semigroup $e^{\sA t}$. 
Thus, the regularity of initial data is naturally measured by means of norms in
intermediate spaces between the domain of the semigroup generator $\sD(\sA)$ and
the state space $Y$; see \eqref{key-2} below.
(Recall from Remark~\ref{r:fractional-powers} that for simplicity of notation we 
write $\sA^\alpha$ instead of $(-\sA)^\alpha$.)
It is this latter formulation, more appropriate in the present context, to require
some adjustment of the arguments used for the proof of Lemma 5.2 in \cite{las-tuff-3}. 
The technical details are the same as the ones in \cite{barbu-etal-1,las-tuff-3}, 
if one takes into account that the solutions considered are generated by the semigroup flow. 
However, for the reader's convenience the proof of Lemma~\ref{l:wave-regularity}
follows below.

\begin{lemma} \label{l:wave-regularity}
Let $(u,w,w_t)$ be the solution to the PDE problem \eqref{e:navierstokes-0} 
corresponding to initial data $y_0 =(u_0,w_0,z_0)\in \sD(\sA^{\alpha})$, 
$0\le \alpha \le 1$, and let $f:=u|_{\Gamma_s}$ so that 
$f\in L_2(0,T;H^{1/2}(\Gamma_s))$.
Then, the $w$ component, which solves the initial/boundary value problem 
\begin{equation}\label{e:open-loop}
\begin{cases}
w_{tt}- {\rm div}\,\sigma(w) =0 & \textrm{in }\; Q_s
\\
\frac{d}{dt} w|_{\Gamma_s}=f & \textrm{on }\; \Sigma_s
\\
w(\cdot,0)=w_0\,,\; w_t(\cdot,0)=z_0\quad  & \textrm{in }\; \Omega_s\,,
\end{cases}
\end{equation}
can be decomposed as $w=w_1+w_2$, where 
$\sigma(w_1)\cdot \nu\in C([0,T],H^{-1/2}(\Gamma_s))$, while
$\sigma(w_2)\cdot \nu\in L_2(0,T;L_2(\Gamma_s))$.
If, in addition, $f\in H^\alpha(\Sigma_s)$, then
$\sigma(w_2)\cdot \nu\in H^\alpha(\Sigma_s)$.
Moreover, the following estimates hold true.
\begin{align}
\|\sigma(w_1)\cdot \nu\|_{C([0,T],H^{-1/2}(\Gamma_s))}^2
& \le C_1 \big(|w_0|_{1,\Omega_s}^2+|z_0|_{0,\Omega_s}^2
+\|f\|_{L_2(0,T;H^{1/2}(\Gamma_s)}\big)
\label{key-1}\\[3pt]
\|\sigma(w_2)\cdot \nu\|_{H^{\alpha}(\Sigma_s)}^2
&\le C_2 \big(\| \sA^{\alpha} y_0\|_H 
+\|f\|_{H^{\alpha}(\Sigma_s)}\big)
\label{key-2}
\end{align}
\end{lemma}

\begin{proof}
The first part of the Lemma follows directly from \cite[Theorem~3.3]{barbu-etal-1}.
Similarly, the regularity inequality \eqref{key-1} is implied by \cite[Lemma 5.2]{las-tuff-3}.
Indeed, this estimate involves the {\em energy level} regularity of the data 
$(w_0,z_0)$, and therefore it is a simple consequence of a more general result---valid 
for any open loop problem like \eqref{e:open-loop}---established in 
\cite[Lemma 5.2]{las-tuff-3}. 
The inequality in \eqref{key-2} involves instead initial data whose regularity is `above' the energy level, and it is thus expressed in terms of the domains of fractional powers of the generator $\sA$ (rather than fractional Sobolev spaces as in \cite{las-tuff-3}). 

In order to establish \eqref{key-2} we appeal, as in \cite{las-tuff-3}, to interpolation  methods. 
Accordingly, it will suffice to prove \eqref{key-2} for $\alpha =0$ and $\alpha = 1$.  
The regularity result corresponding to $\alpha =0$ has been shown in 
\cite[Theorem~3.3]{barbu-etal-1}; thus, it remains to analyze the case $\alpha =1$. 
Note that \eqref{key-2} with $\alpha =1$ provides $H^1(\Sigma_s)$-regularity of the component $\sigma(w_2)$ in terms of $\|\sA y_0\|_H + \|f\|_{H^1(\Sigma_s)}$.
As in \cite{las-tuff-3} and \cite{barbu-etal-1}, we exploit the fact that
$w_2 = \chi w$, where $\chi$ is supported in the microlocal hyperbolic sector;
see \cite[Section~5]{barbu-etal-1}. 
This means that the tangential boundary regularity is implied by the time regularity alone. 
In other words, $H^1(\Sigma_s)$-regularity is equivalent to $H^1(0,T;L_2(\Gamma_s))$ for functions supported in the hyperbolic sector; see \cite[Section~5]{barbu-etal-1} for 
more details. 

In view of the above remarks, we need to show the following estimate:
\begin{equation}\label{key-3}
\Big\|\frac{d}{dt} \sigma (w_2) \cdot \nu \Big\|_{L_2(\Sigma_s)}
\leq C \big(\|\sA y_0 \|_H + \|f\|_{H^1(\Sigma_s)}\big)\,.
\end{equation}
In turn, this inequality expresses the so called ``hidden regularity''
for the wave equation, written in the variable
$\bar{w} = w_t$, where $\bar{w}(0) = z_0, \bar{w}_t(0) = \textrm{div} \sigma (w_0)$,
$\bar{w}|_{\Gamma_s} = f$.
Namely, $\bar{w}$  satisfies 
\begin{equation*}
\begin{cases}
\bar{w}_{tt} - \textrm{div}\,\sigma (\bar{w}) =0 
\nonumber \\
\bar{w}|_{\Gamma_s} = f 
\nonumber \\
\bar{w}(\cdot,0) = z_0\,, \;\bar{w}_t(\cdot,0) = \textrm{div}\,\sigma(w_0)\,.
\end{cases}
\end{equation*}
Owing to Theorem~2.1 in \cite{lasiecka-lions-triggiani}, we obtain 
\begin{equation}\label{new-key}
\|\sigma(\bar{w}) \cdot \nu\|_{L_2(\Sigma_s )} \le C \,\big(\|f\|_{H^1(\Sigma_s)} 
+\|\bar{w}(\cdot,0)\|_{H^1(\Omega)} + \|\bar{w}_t(\cdot,0)\|_{L_2(\Omega)}\big)\,,
\end{equation}
subject to the {\em compatibility condition} $\bar{w}(\cdot,0)|_{\Gamma_s} = f(\cdot,0)$. 
The conditions needed for applicability of the said regularity are 
\begin{equation}\label{e:compatibility}
\begin{cases}
z_0 \in H^1(\Omega)\,,\; \textrm{div}\,\sigma (w_0) \in L_2(\Omega_s)\,, \;
f\in H^1(\Sigma_s)\,, 
\\[2pt]
\textrm{with the compatibility condition $z_0(\cdot)|_{\Gamma_s} = f(\cdot,0)$.}
\end{cases}
\end{equation}
Then, \eqref{new-key} reduces to:
\begin{equation}\label{new-key-1}
\|\sigma(\bar{w}) \cdot \nu \|_{L_2(\Sigma_s)} \le C\,
\big(\|f\|_{H^1(\Sigma_s)} +
\|\textrm{div}\,\sigma(w_0)\|_{L_2(\Omega_s)} + \|z_0\|_{H^1(\Omega_s)}\big)\,.
\end{equation}
Notice now that the conditions \eqref{e:compatibility} are satisfied by any 
$y_0=(u_0,w_0,z_0) \in \sD(\sA)$; indeed, from the very definition \eqref{e:dom-sA}
(of the domain of the generator $\sA$) it follows 
\begin{equation*}
z_0 \in H^1(\Omega)\,,\;
\textrm{div}\,\sigma (w_0)\in L_2(\Omega_s)\,,\;
z_0|_{\Gamma_s} = u(\cdot,0)|_{\Gamma_s} = f(\cdot,0)\,,
\end{equation*}
where the last equality follows from the identification $f = u|_{\Gamma_s}$.
This yields
\begin{equation}\label{new-key-2}
\|\sigma(\bar{w}) \cdot \nu\|_{L_2(\Sigma_s)}
\le C \,\big(\|f\|_{H^1(\Sigma_s)}+\|\sA y_0\|_H\big)\,,
\end{equation}
which in turn implies \eqref{key-3}, after restricting $w$ to $\chi w = w_2$
(we use here the same arguments as in \cite{barbu-etal-1}: $w_2=  \chi w$ satisfies 
the wave equation with microlocal support in the hyperbolic sector). 
Interpolation with the inequality corresponding to $\alpha =0$ finally establishes 
\eqref{key-2}.
\end{proof}

\smallskip 
The sharp results of Lemma~\ref{l:wave-regularity} yield the following improved boundary regularity for the fluid velocity field.

\begin{lemma}
\label{l:key}
Consider the PDE problem~\eqref{e:navierstokes-0}, and take initial data
$y_0 =(u_0,w_0,z_0)\in \sD(\sA^{\alpha})$, $0\le \alpha < \frac14$.
Then, for any $T < \infty $ we have
$u|_{\Gamma_s}\in H^{\alpha}(\Sigma_s)$ and the following estimate holds true:
\begin{equation}\label{e:improved}
\|u\|_{H^{\alpha}(\Sigma_s)}\le C \|\sA^{\alpha} y_0\|_Y\,.
\end{equation}
\end{lemma}
The proof of this Lemma is identical with the one given in \cite[Lemma~5.3]{las-tuff-3}, 
after taking into consideration the estimates derived in Lemma~\ref{l:wave-regularity}.

\section{Proof of the main result}
\label{s:proof}
This section is devoted to the proof of Theorem~\ref{t:main}.
As we shall see, the regularity of boundary traces for the fluid component of the system 
will follow combining the exceptional boundary regularity of the hyperbolic component of the system (Lemma~\ref{l:wave-regularity}), the smoothing properties of the (fluid) analytic semigroup, and the theory of interpolation spaces.

To begin with, let us recall from \cite{las-tuff-3} the basic regularity of the 
Neumann map, which will used repeatedly in the proof of our main result.


\begin{lemma}[\cite{las-tuff-3}] 
\label{l:neumann-trace}
For the map $N:L_2(\Gamma_s)\to H$ defined in \eqref{e:neumann-map} the following 
regularity results holds.
\\
$(i)$ $N\in \sL(L_2(\Gamma_s),\sD(A^{3/4-\delta}))\cap
\sL(H^{-1/2}(\Gamma_s),\sD(A^{1/2}))$ for any $\delta$, $0<\delta<\frac34$.\\ 
$(ii)$ One has $N^*Au = -u|_{\Gamma_s}$, $u\in V$, where the adjoint of $N$ 
is computed with respect to the $L_2$-topology.
\end{lemma}

This Lemma is well known in the context of classical elliptic equations, where
$N$ is the associated Green (Neumann) map. 
The same result has been proved for the abstract operators (defined by \eqref{e:fluid-operator} and \eqref{e:neumann-map}) which arise in the PDE problem under investigation; see \cite[Proposition~4.3]{las-tuff-3}.
For the reader's convenience we recall few details. 
 
The composition $N^*A$ ($N^*$ being the $L_2$-adjoint of $N$)   
is computed from 
\begin{equation*}
\langle N^*A u,g\rangle = (Au, Ng)_f = - \langle u|_{\Gamma_s},g\rangle
\quad\textrm{for all $u \in \sD(A)$, $g \in L_2(\Gamma_s)$.}
\end{equation*} 
Since $u|_{\Gamma_s}$ is bounded on $V = \sD(A^{1/2})$, the above identity can be 
extended by density to all of $V$, even more, to $\sD(A^{1/4 + \delta})$. 
This gives the identification in part (ii) of the Lemma.  

\begin{remark}
{\rm 
In view of part (i) of Lemma \ref{l:neumann-trace}, we can define the 
operator $AN$ as acting from $L_2(\Gamma_s)$ into $[\sD(A^{1/4 + \delta})]'$. 
On the other hand, the analyticity of $e^{At}$ yields the standard bound 
$\|A^{\theta} e^{At}\|_{\sL(H)} \le C\, t^{-\theta}$, $0< t \le 1$. 
This estimate, combined with the invariance of the semigroup $e^{At}$ on the domains of 
fractional powers of $A$, allows to define the convolution  
$\int_0^t e^{A(t-s)} A N g(s) ds$ with values (for each $t$) in
$\sD(A^{3/4 - \delta}) \subset V$, for any $g\in L_{\infty}(0,T;L_2(\Gamma_s))$. 
Thus, $-N^* A \int_0^t e^{A(t-s)} A N g(s) ds$ can be identified---on the strength of 
Lemma~\ref{l:neumann-trace}---with the trace on the boundary of 
$\int_0^t e^{A(t-s)} A N g(s) ds$. 

In the remainder of this section we will use freely the above identifications. 
}
\end{remark}
    
{\em Proof of Theorem~\ref{t:main}. }
Our starting point is the equation satisfied by $u(\cdot)$, namely
$u_t=(A-L)u+AN \sigma(w)\cdot\nu$, whose evolution/integral equivalent is 
\begin{equation} \label{e:mild-1}
u(t) = e^{At}u_0 - \int_0^t e^{A(t-s)} Lu(s) \,ds
+ \int_0^t e^{A(t-s)} AN \sigma(w(s))\cdot\nu \,ds\,.
\end{equation}
We claim that the expected decomposition \eqref{e:decomposition}, along 
with the various regularity properties listed in the statement of 
Theorem~\ref{t:main}, is the one with
\begin{equation} \label{e:terms-decomposition}
u_1(t):= e^{At}u_0\,, \quad 
u_2(t):= -\int_0^t e^{A(t-s)} Lu(s)\,ds+\int_0^t e^{A(t-s)} AN \sigma(w(s))\cdot\nu\,ds\,.
\end{equation}

\noindent
1. Since by Lemma~\ref{l:neumann-trace} we know that $u|_{\Gamma_s}=-N^*Au$,
the splitting \eqref{e:terms-decomposition} yields as well
\begin{align}
& u_1|_{\Gamma_s} =-N^*Au_1(t)= -N^*A e^{At}u_0\,, 
\nonumber\\
& u_2|_{\Gamma_s} =-N^*Au_2(t)=
\underbrace{N^*A\int_0^t e^{A(t-s)} Lu(s)\,ds}_{U_1(t)}
-\underbrace{N^*A\int_0^t e^{A(t-s)} AN \sigma(w(s))\cdot\nu \,ds}_{U_2(t)}\,, 
\label{e:traces-expression}
\end{align}
respectively.
The term $u_1$ yields a pretty straightforward regularity estimate.
Indeed, combining the (analytic) estimates \eqref{e:analytic-estimate} with the 
regularity of the Neumann map recalled in Lemma~\ref{l:neumann-trace}, $(ii)$,
we obtain
\begin{equation}
\begin{split}
\|u_1(t)\|_{L_2(\Gamma_s)}= \|N^*A e^{At}u_0\|_{L_2(\Gamma_s)}
=\|N^*A^{3/4-\delta}A^{1/4+\delta} e^{At}u_0\|_{L_2(\Gamma_s)}
\\[2mm]
\le \|N^*A^{3/4-\delta}\|\, \|A^{1/4+\delta} e^{At}u_0\|_H
\le C\,t^{-1/4-\delta}\|u_0\|\le C_1\,t^{-1/4-\delta}\|y_0\|_Y\,,
\end{split}
\end{equation}
i.e. the ``singular'' estimate \eqref{e:singular-estimate} holds true.
In particular,
\begin{equation} \label{e:result-u_1}
u_1|_{\Gamma_s}\in L_{4-h}(0,T;L_2(\Gamma_s)) \quad
\textrm{with arbitrarily small $h>0$.}
\end{equation}
Thus, to establish \eqref{e:u-reg}, it remains to demonstrate that the boundary trace of the second component $u_2$ on $\Gamma_s$ is well defined and satisfies 
\eqref{e:u_2-reg}.
Taking account of \eqref{e:traces-expression}, since the regularity of $U_2$ has
been dealt with in \cite{bucci-las-cvpde}, we will be done once we show
\begin{equation} \label{e:U_1-reg}
U_1\in L_p(0,T;L_2(\Gamma_s)) \qquad \forall p\ge 1\,.
\end{equation}
We consider the convolution integral $z(t) = \int_0^t e^{A(t-s)} Lu(s)\,ds$,
and recall from \cite[Theorem~3.2]{barbu-etal-1} the additional property 
$u\in L_2(0,T;V)$, which implies $\sL u \in L_2(0,T;L_2(\Omega_f))$ and  
$Lu = P_H \sL u \in L_2(0,T;H)$.
Thus, by parabolic regularity theory we deduce $z\in L_2(0,T;\sD(A))\cap H^1(0,T;H)$;
therefore, by interpolation, 
\begin{equation*}
z\in H^\alpha(0,T;\sD(A^{1-\alpha})) \qquad \forall \alpha\in (0,1)\,.
\end{equation*}
In particular, we may take $\alpha = \frac34-\delta$ to find
\begin{equation*} 
z\in H^{3/4-\delta}(0,T;\sD(A^{1/4+\delta}))\,,
\end{equation*}
and since 
\begin{equation*}
U_1(t)=[N^*A^{3/4-\delta}]\, A^{1/4+\delta}\,z(t)\,,
\end{equation*}
we obtain 
\begin{equation} \label{e:U_1-again}
U_1\in H^{3/4-\delta}(0,T;L_2(\Gamma_s))\subset C([0,T],L_2(\Gamma_s))\,,
\end{equation}
where the inclusion follows by Sobolev embedding theorems in one dimension.
Thus, the validity of \eqref{e:U_1-reg} is ensured by \eqref{e:U_1-again}.

\smallskip
The analysis of the summand $U_2(t)$ critically relies on the sharp regularity theory 
for the elastic component of the system, as given by the first part of 
Lemma~\ref{l:wave-regularity}. Specifically, on the basis of 
\begin{equation*} 
\sigma(w)\cdot\nu \in C([0,T],H^{-1/2}(\Gamma_s))+ L_2(0,T;L_2(\Gamma_s))\,,
\end{equation*}
it was shown in \cite{bucci-las-cvpde} that 
\begin{equation} \label{e:result-U_2}
U_2\in C([0,T],L_2(\Gamma_s)) +  L_p(0,T;L_2(\Gamma_s)) 
\subset L_p(0,T;L_2(\Gamma_s)) \qquad \forall p\ge 1\,.
\end{equation}
Consequently, combining \eqref{e:result-U_2} with \eqref{e:U_1-reg} we obtain
\eqref{e:u_2-reg}, as required. 
This concludes the proof of all assertions in item 1.

\medskip
\noindent 
2. Let now $y_0\in \sD({\sA}^\epsilon)$, $\epsilon>0$.
Aiming to show the validity of \eqref{e:u_2-reg-improved}, we recall 
$u_2|_{\Gamma_s}=-N^*Au_2=U_1-U_2$, with $U_i$ defined in 
\eqref{e:traces-expression}.
We have seen in step 1  
that $U_1\in C([0,T],L_2(\Gamma_s))$ when $y_0$ 
just belongs to $Y$; the same holds true, {\em a fortiori}, 
in the present case.

As for $U_2$, its regularity analysis when $y_0\in \sD({\sA}^\epsilon)$ has been pursued 
in \cite{bucci-las-cvpde}. We just note that since 
\begin{equation*}
y_0\in \sD({\sA}^\epsilon)\,,
\end{equation*}
according to Lemma~\ref{l:key} $u\in H^\epsilon(\Sigma_s)$, provided 
that $\epsilon<\frac14$.
Therefore, we may apply the second part of Lemma~\ref{l:wave-regularity}
($f=u|_{\Gamma_s}$), thereby obtaining
\begin{equation}\label{e:additional}
\sigma(w)\cdot \nu\in C([0,T],H^{-1/2}(\Gamma_s)) + H^{\epsilon}(\Sigma_s)\,.
\end{equation}
Owing to \eqref{e:additional}, a decomposition of $U_2$ is necessitated, 
say $U_2=U_{21}+U_{22}$.
Careful calculations made in \cite{bucci-las-cvpde} showed that 
now both terms $U_{2i}$ belong to $C([0,T],L_2(\Gamma_s))$, $i=1,2$, which combined 
with \eqref{e:U_1-again} finally yields the conclusion \eqref{e:u_2-reg-improved}.

\medskip
\noindent
3. We turn our attention to the time derivative $u_t$.
In order to prove \eqref{e:u_t-reg}, we compute from \eqref{e:mild-1} 
(in the dual space $[\sD(A)]'$) 
\begin{equation}
u_t(t) := v_1(t)+v_2(t)= Ae^{At}u_0 + v_2(t)\,,
\end{equation}
where, initially, $y_0\in Y$ and $v_2$ reads as 
\begin{align*} 
& v_2(t)=-\big[Lu(t) + A\int_0^t e^{A(t-s)} Lu(s)\,ds\Big]
\\
& \quad + \Big[A\int_0^t e^{A(t-s)} AN \sigma(w(s))\cdot\nu\,ds 
+ AN \sigma(w(t))\cdot\nu\Big]=:-V_1(t)+V_2(t)\,.
\end{align*}
As $V_1$ is nothing but $z_t$ ($z$ is defined in step 1), and according to parabolic regularity $z\in H^1(0,T;H)$, we immediately obtain
\begin{equation}\label{e:V_1}
V_1\in L_2(0,T;H)\subset L_2(0,T;L_2(\Omega_f))\,.
\end{equation}
The more challenging analysis of $V_2$ yields, utilizing once again 
Lemma~\ref{l:wave-regularity},
\begin{equation}\label{e:V_2}
V_2\in L_2(0,T;[D(A^{1/2+\delta_1})]')
= L_2(0,T;[H^{1+2\delta_1}(\Omega_f)]')\,,
\end{equation}
for any $0<\delta_1<\frac14$; see \cite{bucci-las-cvpde} for all details.
Combining \eqref{e:V_1}--\eqref{e:V_2} leads us to the estimate 
\begin{equation}\label{e:v_2-first}
y_0\in Y \Longrightarrow v_2\in L_2(0,T;[H^{1+2\delta_1}(\Omega_f)]')\,,
\qquad 0<\delta_1<\frac14\,.
\end{equation}

\smallskip
When $y_0\in \sD(\sA)$, we rewrite $v_2$ in a different fashion, namely
\begin{multline}\label{e:decomp-2}
v_2(t)=-\Big[e^{At}Lu_0 + \int_0^t e^{A(t-s)}Lu_s(s)\,ds\Big]
\\
+ \Big[\int_0^t e^{A(t-s)} AN \sigma(w_s(s))\cdot\nu\,ds
+Ae^{At}N \sigma(w(0))\cdot\nu\Big]=:-Z_1(t)+Z_2(t)\,.
\end{multline}
According to well-posedness theory in \cite[Theorem~3.2]{barbu-etal-1}
$y_0\in \sD(\sA)$ implies now $y_t=(u_t,w_t,w_{tt})\in C([0,T];Y)$, with 
the additional regularity $u_t\in L_2(0,T;V)$. 
Consequently, $Lu_t\in L_2(0,T;H)$ and 
\begin{equation} \label{e:Z_12}
\int_0^{\cdot} e^{A(\cdot-s)}Lu_s(s)\,ds\in L_2(0,T;\sD(A))
\subset L_2(0,T;H^2(\Omega_f))\,.
\end{equation} 
On the other hand, with $y_0\in \sD(\sA)$ we have $u_0\in V$ so that 
$Lu_0\in H$ and in particular
\begin{equation*} 
e^{A\cdot}Lu_0\in L_2(0,T;\sD(A^{1/2-\delta}))= L_2(0,T;H^{1-2\delta}(\Omega_f))\,, 
\end{equation*}
which is weaker than \eqref{e:Z_12}; consequently, 
\begin{equation} \label{e:Z_1}
Z_1\in L_2(0,T;H^{1-2\delta}(\Omega_f))\,, \qquad 0<\delta<\frac12\,.
\end{equation}

\smallskip
Consider next the term $Z_2$, and notice that $y_0\in \sD(\sA)$ implies,
in particular, $(w_0,z_0)\in H^2(\Omega_s)\times H^1(\Omega_s)$.
By standard semigroup arguments and according to Lemma~\ref{l:wave-regularity}
we know that
\begin{equation*}
\sigma(w_t)\cdot \nu \in C([0,T],H^{-1/2}(\Gamma_s))\,+\,L_2(0,T;L_2(\Gamma_s))\,.
\end{equation*}
In light of the above result, it was shown in \cite{bucci-las-cvpde} that 
the convolution integral
\begin{equation*}
Z_{21}:=\int_0^{\cdot} e^{A(\cdot-s)} AN \sigma(w_s(s))\cdot\nu\,ds
\end{equation*}
satisfies 
\begin{equation*}
Z_{21}\in C([0,T],D(A^{1/2-\delta}))\,+\,L_{4}(0,T;D(A^{1/2-\delta}))\,,
\qquad 0<\delta<\frac12\,,
\end{equation*}
which gives 
\begin{equation}\label{e:Z_21}
Z_{21}\in L_4(0,T;D(A^{1/2-\delta}))\,, \qquad 0<\delta<\frac12\,.
\end{equation}

On the other hand, still with $y_0\in D(\sA)$, one has just 
$\sigma(w)\cdot\nu\in C([0,T],H^{-1/2}(\Gamma_s))$. 
Rewrite the summand $Z_{22}(t):= Ae^{At}N \sigma(w(0))\cdot\nu$
in \eqref{e:decomp-2} as
\begin{equation*}
Z_{22}(t)=A^{1/2}e^{At}\,\big[A^{1/2}N \sigma(w(0))\cdot\nu\big]\,,
\end{equation*}
and use once more Lemma~\ref{l:neumann-trace}, along with the 
singular estimates pertaining to analytic semigroups, to find
\begin{equation}\label{e:Z_22}
Z_{22}\in L_q(0,T;D(A^{1/2-\delta}))\,, \quad \textrm{provided that $q(1-\delta)<1$.}
\end{equation}
for any $\delta<\frac12$. 
Recalling $v_2=-Z_1+Z_{21}+Z_{22}$, with the requirements 
\eqref{e:Z_1}, \eqref{e:Z_21} and \eqref{e:Z_22}, we finally obtain
\begin{equation}\label{e:v_2-second}
y_0\in \sD(\sA)\Longrightarrow v_2\in L_{q_2}(0,T;D(A^{1/2-\delta_2})) 
= L_{q_2}(0,T;H^{1-2\delta_2}(\Omega_f))
\end{equation}
for any $0<\delta_2<\frac12$, where $q_2\in (1,2)$ depends on $\delta_2$; 
more precisely, 
\begin{equation}\label{e:sobolev-exponent_2}
q_2< \frac1{1-\delta_2}\,. 
\end{equation}

\smallskip
We proceed now as in the conclusion of step $(iii)$ in 
\cite[Proof of Theorem~2.9]{bucci-las-cvpde}.
Combining \eqref{e:v_2-second} with \eqref{e:v_2-first}, gives by interpolation
\begin{equation}\label{e:interpolation}
y_0\in D(\sA^{1-\theta})\Longrightarrow  v_2\in L_{q_2}(0,T;W)\,, 
\end{equation}
where $q_2$ is as in \eqref{e:sobolev-exponent_2} and $W$ is the interpolation space 
\begin{equation*}
W= (H^{1-2\delta_2}(\Omega_f),[H^{1+2\delta_1}(\Omega_f)]')_{\theta}
\equiv H^{s}(\Omega_f)\,, 
\end{equation*}
if
\begin{equation*}
s=(1-\theta)(1-2\delta_2)-\theta(1+2\delta_1)=
1-2\delta_2-2\theta(1+\delta_1-\delta_2)\ge 0\,;
\end{equation*}
see \cite[Chapter~1, Theorem~12.5]{lions-magenes}.
Notice that by taking, for instance, $\delta_1=\delta_2=:\delta$, 
one has $s\ge 1/2$ provided that
\begin{equation}\label{e:theta-constraint}
\theta+\delta \le \frac14\,.
\end{equation}
In this case $v_2$ takes values in $H^s(\Omega_f)$ with $s\ge 1/2$ and hence its 
trace on $\Gamma_s$ is well defined. 
Notice that, in view of the constraint \eqref{e:theta-constraint},
we need to require $0<\theta<\frac14$.
Consequently, given any $\theta$ such that $0<\theta<\frac14$, 
choosing, e.g., $\delta= 1/4 -\theta$ in view of \eqref{e:theta-constraint}, 
from \eqref{e:interpolation} it follows
\begin{equation}\label{e:trace_1}
y_0\in D(\sA^{1-\theta})) \Longrightarrow
N^*Av_2\in L_{q_2}(0,T;L_2(\Gamma_s)) \qquad \forall q_2< \frac{4}{3+4\theta}\,.
\end{equation}

It remains to establish the regularity of the first summand
$N^*Av_1(t)= N^*Ae^{At}Au_0$ when $y_0\in D(\sA^{1-\theta})$.
In this case $u_0\in (H^1(\Omega_f),L_2(\Omega_f))_\theta=H^{1-\theta}(\Omega_f)$, 
and it is not difficult to show 
\begin{equation}\label{e:trace_2}
y_0\in D(\sA^{1-\theta}) \Longrightarrow
N^*Av_1\in L_{q_1}(0,T;L_2(\Gamma_s)) \qquad \forall q_1< \frac{4}{3+2\theta+4\epsilon}
\end{equation}
where $\epsilon$ is independent of $\theta$ and can be taken arbitrarily small.
(Notice that the Sobolev exponent $q_1$ belongs to $(1,2)$, as well.)
Thus, in view of \eqref{e:trace_1} and \eqref{e:trace_2} we finally obtain 
\begin{equation}\label{e:trace}
y_0\in D(\sA^{1-\theta}) \Longrightarrow
u_t|_{\Gamma_s} \in L_{q}(0,T;L_2(\Gamma_s)) 
\qquad \forall q< \frac{4}{3+4\theta}\,,
\end{equation}
which concludes the proof.
\qed

\section{Application to the quadratic optimal control problem}
\label{s:application}
In this section we briefly discuss the implications of the obtained trace regularity 
results in the study of the optimal control problem (with quadratic functionals)
associated with the fluid-solid interaction model \eqref{e:navierstokes-0}.
Thus, we consider the initial/boundary value problem
\begin{equation} \label{e:navierstokes-1}
\left\{ \hspace{1mm}
\begin{split} 
& u_t-{\rm div}\,\epsilon(u) + \sL u
+\nabla p= 0 & &\textrm{in }\; Q_f:= \Omega_f\times (0,T)
\\
& {\rm div}\, u=0 & & \textrm{in }\; Q_f
\\
& w_{tt} - {\rm div}\,\sigma(w)=0 & &\textrm{in }\; Q_s:= \Omega_s\times (0,T)
\\
& u=0 & & \textrm{on }\; \Sigma_f:= \Gamma_f\times (0,T)
\\
& w_t=u & & \textrm{on }\; \Sigma_s:= \Gamma_s\times (0,T)
\\
& \sigma(w)\cdot \nu=\epsilon(u)\cdot\nu - p\nu-g & & 
\textrm{on }\; \Sigma_s
\\
& u(0,\cdot)=u_0 & & \textrm{in }\; \Omega_f
\\
& w(0,\cdot)=w_0\,, \quad w_t(0,\cdot)=z_0 & & \textrm{in }\; \Omega_s\,,
\end{split}
\right.
\end{equation}
which includes a function $g$ acting as a control on the interface $\Gamma_s$.
The PDE problem \eqref{e:navierstokes-1} readily corresponds to an abstract (linear) control system, that is
\begin{equation} \label{e:control-system}
\begin{cases}
y'= \sA y + \sB g & {\rm in } \quad [\sD({\sA}^*)]'
\\
y(0)=y_0\in Y
\end{cases}\,,
\end{equation}
where the (free dynamics) operator $\sA$ is the one in \eqref{e:sA}, while the  (control) operator $\sB:L_2(\Gamma_s)=:U\to [\sD(\sA)]'$ is given by 
\begin{equation}\label{e:sB}
\begin{array}{l}
\sB=
\left(
\begin{array}{c}
AN\\
0\\
0
\end{array}
\right)\,. 
\end{array}
\end{equation}
Semigroup well-posedness of the uncontrolled problem corresponding to 
\eqref{e:control-system} 
still follows from \cite[Proposition~3.1]{barbu-etal-1}. 
The membership $\sB\in \sL(U,[\sD(\sA)]')$ is a consequence of $AN$ being bounded from $L_2(\Gamma_s)$ into $V'$ (yet unbounded from $L_2(\Gamma_s)$ into $H$); 
see \cite{las-tuff-3} for more details.

We associate to \eqref{e:control-system} a quadratic functional
\begin{equation} \label{e:functional}
J(g) = \int_0^T\big(|\sR y(t)|_Z^2 + |g(t)|_U^2\big)\,dt
\end{equation}
where $g$ varies in $L_2(0,T;U)$, and initially $\sR\in \sL(Y,Z)$ 
($\sR$ is called {\em observation} operator, $Z$ is the observations space).
Then, the (optimal) quadratic control problem over a finite horizon is stated
as follows.

\begin{problem}[The optimal control problem]
Given $y_0\in Y$, seek a control function $g\in L_2(0,T;U)$ which minimizes the 
cost functional \eqref{e:functional}, where $y(t)=y(t;y_0,g)$ is the solution to 
\eqref{e:control-system} corresponding to $g$.
\end{problem}
As in the classical Linear Quadratic (LQ) problem in a finite dimensional
context, solvability of the optimal control problem defined above does not just
mean existence of a (unique) minimizer $g^*$.
We aim to obtain the {\em feedback} representation of the optimal control, defined for any $y_0$ in the state space $Y$,
along with well-posedness of the corresponding operator Riccati equations.
On the other hand, it is well known that in the case of boundary control systems, 
such as the present \eqref{e:control-system}, a major difficulty comes from the fact 
that the {\em gain} operator $\sB^* P(t)$ which occurrs in the feedback formula 
\begin{equation*}
g^*(t) = -\sB^* P(t) y^*(t) \qquad \textrm{a.e. in $[0,T]$,}
\end{equation*}
may not be defined, or may be {\em unbounded}.

Let us recall that owing to the well-known representation of the Riccati operator 
$P(t)$ (involving the {\em evolution map} $\Phi(\cdot,\cdot)$ associated with the optimal dynamics)
\begin{equation*}
P(t)z = \int_t^T e^{\sA^*(s-t)} R^*R \Phi(s,t)z\, ds\,, \quad z\in Y\,,
\end{equation*} 
it is immediately seen that boundedness of the operator $\sB^* P(t)$ 
will crucially depend upon the regularity properties of the operator
\begin{equation*}
\sB^*e^{\sA^*t}:\sD(\sA^*)\to U=L_2(\Gamma_s)\,.
\end{equation*}
Thus, a PDE intepretation of the composition $\sB^*e^{\sA^*t}$ is called for.
This is obtained---for a given control system \eqref{e:control-system}---by computing 
$\sB^*$ and by noting that $e^{\sA^*t}y_0$ solves the {\em uncontrolled} Cauchy problem 
$\hat{y}'=\sA^* \hat{y}$, $\hat{y}(0)=y_0$.
These calculations eventually yield the distinct {\em trace} operator corresponding 
to $\sB^*e^{\sA^*t}$.
(This connection appeared clear since the former studies of the LQ-problem for simple hyperbolic equations; see, e.g., \cite{las-trig-books_2} and \cite{bddm}.)
In the present case we obtain, in particular,
\begin{equation}
\sB^*e^{\sA^*\cdot}y_0 = [N^*A,0,0]\hat{y} =N^*A\hat{u}=-\hat{u}|_{\Gamma_s}\,,
\end{equation}
where $\hat{y}=(\hat{u},\hat{w},\hat{w}_t)$ solves the initial/boundary value problem corresponding to the afore-said adjoint system, which is almost identical with the PDE problem \eqref{e:navierstokes-0}.
This explains why a study of the boundary regularity of the fluid velocity 
field on the interface is necessitated.

A situation where one can conclude that $\sB^* P(t)$ is bounded on the state space $Y$ 
is when the pair $(\sA,\sB)$ yields a `singular' estimate for the operator 
$e^{\sA t}\sB$ near $t=0$; see \cite{las-cbms}, \cite{las-trig-se-2}, and \cite{las-tuff-2}.
It was shown in \cite{las-tuff-3} that the PDE problem \eqref{e:navierstokes-1} 
does give rise to a singular estimate; however, this holds true provided the observation operator $\sR$ possesses a suitable smoothing property.

In contrast, in light of the regularity results established in Theorem~\ref{t:main}, 
we are able to show that the abstract control system \eqref{e:control-system}
corresponding to the initial/boundary value problem \eqref{e:navierstokes-1} 
falls in the more general class of systems introduced in \cite{abl-2}. 
Indeed, all the regularity results proved in Theorem~\ref{t:main}---which 
yield relative assertions for a similar (dual) PDE problem---are the PDE 
counterpart of the control-theoretic properties (Hypotheses~2.1--2.2 therein) 
required to apply the theory of \cite{abl-2}.
We do not discuss the proof, which follows the lines of 
\cite[Proof of Theorem~2.6]{bucci-las-cvpde}, but let us remark that
while the regularity of the boundary traces of $u$ on $\Gamma_s$ corresponds to 
the regularity of the operator $\sB^*e^{\sA^*t}$, the obtained regularity for 
$u_t$ in \eqref{e:u_t-reg} yields a sought-after and most challenging property for 
$\sB^*e^{\sA^*t}{\sA^*}^\theta$ (that is condition (3.23) in \cite{abl-2}, or
(4.4) in \cite{bucci-las-cvpde}).
%
We just note that the latter PDE interpretation follows from the following key 
observation:
\begin{equation*}
\sB^*e^{\sA^*t}{\sA^*}^\theta y_0 = 
\sB^*\frac{d}{dt} e^{\sA^*t} ({\sA^*}^{\theta-1} y_0) =
N^*A\hat{u}_t=-\hat{u}_t|_{\Gamma_s}\,,
\end{equation*}
where $(\hat{u},\hat{w},\hat{w}_t)$ solves the same (adjoint) PDE problem as before, yet with initial state $\xi_0\equiv {\sA^*}^{\theta-1} y_0\in \sD({\sA^*}^{1-\theta})$ 
(while time derivatives are understood distributionally); 
the reader is referred to \cite[Section~4.1]{bucci-las-cvpde} for more details.

Thus, in order for all the hypotheses needed to apply \cite[Theorem~2.3]{abl-2}
to be fulfilled, we just need to assume the one which pertains to the observation 
operator $\sR$; namely,
\begin{equation} \label{h:observation}
{\sR}^*\sR\in \sL(\sD({\sA}^{\epsilon}),\sD({{\sA}^*}^{\epsilon})) 
\end{equation}
for some $\epsilon\in (0,1/4)$.
Then, under \eqref{h:observation} all the conclusions of \cite[Theorem~2.3]{abl-2} follow. 
In particular, one obtains the feedback synthesis of the optimal control for all 
$y_0$ in the finite energy space $Y$, with a gain operator which is 
{\em densely defined}, that is enough to get well-posedness of the 
corresponding differential Riccati equations.

\begin{remark}\label{r:nonsmoothing-observation}
\begin{rm}
We finally note that while a case of particular interest is $\sR=I$,
the requirement \eqref{h:observation} does not `force' smoothing effects of 
the observation operator. 
Indeed, in the present case $\sD({\sA^*}^{\epsilon})$ coincides with
$\sD(\sA^{\epsilon})$ for $\epsilon$ sufficiently small, and hence 
\eqref{h:observation} is satisfied by any operator which `maintains' 
regularity, such as the identity $\sR=I$;  
see \cite[Remark~2.7]{bucci-las-cvpde} for a more detailed explanation.
This answers the final question raised in \cite[Remark~6.1]{las-tuff-3}.
\end{rm}
\end{remark}

\section*{Acknowledgements}
We wish to thank the anonymous referees for the careful reading and the useful comments
on the former version of the manuscript.
We are specially grateful to a referee, whose remarks led us to amend the proof 
of Theorem~\ref{t:main} (step~2), giving rise to novel additions in section~\ref{s:known-results}.


\end{document}